\documentclass[12pt,a4paper,english,intoc,bibliography=totoc,index=totoc,BCOR10mm,captions=tableheading,titlepage,fleqn]{scrbook}
\usepackage{lmodern}

\usepackage[T1]{fontenc}
\usepackage[latin9]{inputenc}
\usepackage{babel}
\usepackage{amsmath}
\usepackage{amssymb}
\usepackage[unicode=true,
 bookmarks=true,bookmarksnumbered=true,bookmarksopen=true,bookmarksopenlevel=2,
 breaklinks=false,pdfborder={0 0 0},backref=false,colorlinks=false]
 {hyperref}
\hypersetup{pdftitle={Abstract},
 pdfauthor={Uwe Stöhr},
 pdfpagelayout=OneColumn, pdfnewwindow=true, pdfstartview=XYZ, plainpages=false}
\usepackage{breakurl}

\makeatletter

\usepackage{calc}

\makeatother

\begin{document}
\noindent \begin{center}
\textbf{\Large On the non-existence of $\kappa$-mad families}
\par\end{center}{\Large \par}

\noindent \begin{center}
{\large Haim Horowitz and Saharon Shelah}
\par\end{center}{\large \par}

\noindent \begin{center}
\textbf{\small Abstract}
\par\end{center}{\small \par}

\noindent \begin{center}
{\small Starting from a model with a Laver-indestructible supercompact
cardinal $\kappa$, we construct a model of $ZF+DC_{\kappa}$ where
there are no $\kappa$-mad families.}%
\footnote{{\small Date: June 12, 2019}{\small \par}



Publication 1168 of the second author%
}
\par\end{center}{\small \par}

\textbf{\large Introduction}{\large \par}

The study of the definability and possible non-existence of mad families
has a long tradition, originating with the paper {[}Ma{]} of Mathias
where it was proven that mad families can't be analytic and that there
are no mad families in the Solovay model constructed from a Mahlo
cardinal. It was later shown by Toernquist that an inaccessible cardinal
suffices for the consistency of this statement ({[}To{]}), and it
was then shown by the authors that the non-existence of mad families
(in $ZF+DC$) is actually equiconsistent with $ZFC$ ({[}HwSh:1090{]}).

The current paper can be seen as a continuation of the line of investigation
of {[}HwSh:1090{]}, as well as of {[}HwSh:1145{]}, where the definability
of $\kappa$-mad families was considered. Recall the following definition:

\textbf{Definition 1: }Let $\kappa$ be an infinite regular cardinal.
A family $\mathcal A \subseteq [\kappa]^{\kappa}$ is $\kappa$-almost
disjoint if $|A\cap B|<\kappa$ for every $A\neq B \in \mathcal A$.
$\mathcal A$ will be called $\kappa$-maximal almost disjoint ($\kappa$-mad)
if $\mathcal A$ is $\kappa$-almost disjoint and can't be extended
to a larger $\kappa$-almost disjoint family.

Assuming the existence of a Laver-indestructible supercompact cardinal
$\kappa$, we constructed in {[}HwSh:1145{]} a generic extension where
$\kappa$ remained supercompact and there are no $\Sigma^1_1(\kappa)-\kappa-$mad
families, thus obtaining a higher analog of Mathias' result.

Our current main goal is to obtain a higher analog of the main result
of {[}HwSh:1090{]}, i.e. for an uncountable cardinal $\theta>\aleph_0$,
we would like to construct a model of $ZF+DC_{\theta}$ where there
are no $\theta$-mad families. As opposed to {[}HwSh:1090{]}, we only
achieve this goal assuming the existence of a supercompact cardinal.
The main result of the paper is the following:

\textbf{Theorem 2: }a. Suppose that $\aleph_0<cf(\theta)=\theta<cf(\kappa)=\kappa \leq \lambda=\lambda^{<\kappa}$
and $\theta$ is a Laver indestructible supercompact cardinal, then
there is a model of $ZF+DC_{<\kappa}+"$there exist no $\theta$-mad
families$"$.

b. If we start from a universe $V$, then the final model $V_1$ will
have the same cardinals and same $H(\theta)$ as $V$.

We shall force with a partial order $\mathbb P$ where the conditions
themselves are forcing notions (this is somewhat similar to {[}Sh:218{]},
{[}HwSh:1093{]} and {[}HwSh:1113{]}, as well as to the recent work
of Viale in {[}Vi{]}, where a similar approach is applied to the study
of generic absoluteness). Forcing with $\mathbb P$ will generically
introduce the forcing notion $\mathbb Q$ that will give us the desired
results. More specifically, we shall fix a Laver-indestructible supercompact
cardinal $\theta$. The conditions in $\mathbb P$ will be elements
from a suitable $H(\lambda^+)$ that are $(<\theta)$-support iterations
along wellfounded partial orders of $(<\theta)$-directed closed forcing
notions satisfying a strong version of $\theta^+$-cc. Given $\bold{q}_1,\bold{q}_2 \in \mathbb P$,
we will have $\bold{q}_1 \leq_{\mathbb P} \bold{q}_2$ when the iteration
given by $\bold{q}_1$ is an {}``initial segment'' (in an adequate
sense) of the iteration given by $\bold{q}_2$. Forcing with $\mathbb P$
will introduce a generic iteration $\bold{q}_G$ given by the union
of $\bold q \in \mathbb P$ that belong to the generic set. In the
further generic extension given by $\bold{q}_G$, we shall consider
$V_1=HOD(\mathcal{P}(\theta)^{<\kappa} \cup V)$ (for an adequate
fixed $\kappa$). We shall then prove that there are no $\theta$-mad
families in $V_1$. In order to prove this fact, we shall consider
towards contradiction a condition $(\bold{q}_0,\underset{\sim}{p_0})$
that forces a counterexample $\mathcal A$, where $\bold{q}_0$ will
be {}``sufficiently closed''. The filter that's dual to the ideal
generated by $\mathcal A$ will then be extended to a $\theta$-complete
ultrafilter (using the Laver-indestructibility of $\theta$), and
we shall obtain a contradiction with the help of an amalgamation argument
over $\bold{q}_0$ using a higher analog of Mathias forcing relative
to this ultrafilter.

The rest of the paper will be devoted to the proof of Theorem 2.

\textbf{\large Proof of the main result}{\large \par}

\textbf{Definition 3: }A. Let $K$ be the class of $\bold q$ that
consist of the following objects with the following properties:

a. $U=U_{\bold q}$ a well-founded partial order whose elements are
ordinals. We let $U^+=U\cup \{ \infty\}$ where $\infty$ is a new
element above all elements from $U$, and for $\alpha 
\in U^+$, we let $U_{<\alpha}=\{ \beta \in U : \beta<_{U} \alpha\}$.

b. An iteration $(\mathbb{P}_{\bold q, \alpha}, \underset{\sim}{\mathbb{Q}_{\bold q, \beta}} : \alpha \in U^+, \beta \in U)=(\mathbb{P}_{\alpha}, \underset{\sim}{\mathbb{Q}_{\beta}} : \alpha \in U^+, \beta \in U)$.
We shall often denote the iteration itself by $\bold q$.

c. $\bold q$ is a $(<\theta)$-support iteration, and in addition:

$(\alpha)$ Each $\underset{\sim}{\mathbb{Q}_{\beta}}$ is a $\mathbb{P}_{\beta}$-name
of a forcing notion whose set of elements is an object $X_{\beta}$
from $V$.

$(\beta)$ Given $\alpha \in U^+$, $p\in \mathbb{P}_{\alpha}$ iff
$p$ is a function with domain $dom(p) \in [U_{<\alpha}]^{<\theta}$
such that $p(\beta)$ is a canonical $\mathbb{P}_{\beta}$-name for
every $\beta \in dom(p)$.

$(\gamma)$ $\leq_{\mathbb{P}_{\alpha}}$ is defined as usual.

$(\delta)$ If $w\subseteq U$ is downward closed (i.e. $\alpha<_U \beta \in w \rightarrow \alpha \in w$)
and $\mathbb{P}_{\bold q, w}=\mathbb{P}_w=\mathbb{P}_{\infty} \restriction w=\{ p\in \mathbb{P}_{\infty} : dom(p) \subseteq w \}$,
then $\mathbb{P}_w \lessdot \mathbb{P}_{\infty}$.

d. In $V^{\mathbb{P}_{\beta}}$, $\underset{\sim}{\mathbb{Q}_{\beta}}$
satisfies $*_{\theta}^{\epsilon}$ for a fixed limit $\epsilon<\theta$,
namely, if $\{p_{\alpha} : \alpha<\theta^+ \} \subseteq \underset{\sim}{\mathbb{Q}_{\beta}}$,
then there is some club $E\subseteq \theta^+$ and a pressing down
function $f: E \rightarrow \theta^+$ such that if $\delta_1,\delta_2 \in E$,
$cf(\delta_1)=cf(\delta_2)$ and $f(\delta_1)=f(\delta_2)$, then
$p_{\delta_1}$ and $p_{\delta_2}$ have a common least upper bound.

e. For $\beta \in U$, the following holds in $V^{\mathbb{P}_{\beta}}$:
If $I$ is a directed partial order of cardinality $<\theta$ and
$(p_s : s\in I) \in \mathbb{Q}_{\beta}^I$ is $\leq_{\mathbb{Q}_{\beta}}$-increasing,
then $\{p_s : s\in I\}$ has a $\leq_{\mathbb{Q}_{\beta}}$-least
upper bound.

B. Let $\leq_K$ be the following partial order on $K$:

$\bold{q}_1 \leq_K \bold{q}_2$ iff the following conditions hold:

a. $U_{\bold{q}_1} \subseteq U_{\bold{q}_2}$ as partial orders.

b. If $U_{\bold{q}_2} \models \alpha<\beta$ and $\beta \in U_{\bold{q}_1}$,
then $\alpha \in U_{\bold{q}_1}$.

c. If $w\subseteq U_{\bold{q}_1}$ is downward closed, then $\mathbb{P}_{\bold{q}_1,w}=\mathbb{P}_{\bold{q}_2,w}$.

d. If $\alpha \in U_{\bold{q}_1}$, then $\underset{\sim}{\mathbb{Q}_{\bold{q}_1,\alpha}}=\underset{\sim}{\mathbb{Q}_{\bold{q}_2,\alpha}}$
(this is well-defined recalling clause (b)).

C. Let $K_{wf}$ be the class of $U$ as in (A)(a), and let $\leq_{wf}$
be the partial order on $K_{wf}$ defined as in clauses (B)(a) and
(B)(b).

We shall now observe some easy basic properties of the objects defined
above:

\textbf{Observation 4: }a. If $(U_{\alpha} : \alpha<\delta)$ is $\leq_{wf}$-increasing,
then $\underset{\alpha<\delta}{\cup}U_{\alpha}$ is a $\leq_{wf}$-least
upper bound for $(U_{\alpha} : \alpha<\delta)$.

b. $\leq_K$ is a partial order on $K$.

c. If $\bold{q}_2 \in K$ and $U_1 \subseteq U_{\bold{q}_2}$ is downward
closed, then there is a unique $\bold{q}_1 \in K$ such that $\bold{q}_1 \leq_K \bold{q}_2$
and $U_{\bold{q}_1}=U_1$.

d. If $(\bold{q}_{\alpha} : \alpha<\delta)$ is $\leq_K$-increasing,
then there is a unique $\bold{q}_{\delta} \in K$ such that $\alpha<\delta \rightarrow \bold{q}_{\alpha} \leq_K \bold{q}_{\delta}$
and $U_{\bold{q}_{\delta}}=\underset{\alpha<\delta}{\cup}U_{\bold{q}_{\alpha}}$.

e. If $U_0, U_1, U_2 \in K_{wf}$, $U_0=U_1 \cap U_2$ and $U_0 \leq_{wf} U_l$
$(l=1,2)$, then there is a unique $U\in K_{wf}$ such that $\underset{l=1,2}{\wedge}U_l \leq_{wf} U$,
$\alpha \in U$ iff $\alpha \in U_1 \vee \alpha \in U_2$ and $\leq_U= \leq_{U_1} \cup \leq_{U_2}$.
We denote this $U$ by $U_1+_{U_0} U_2$.

f. If $\bold{q}_0, \bold{q}_1, \bold{q}_2 \in K$, $\bold{q}_0 \leq_K \bold{q}_l$
$(l=1,2)$ and $U_{\bold{q}_0}=U_{\bold{q}_1} \cap U_{\bold{q}_2}$,
then there is a unique $\bold q \in K$ such that $\underset{l=1,2}{\wedge} \bold{q}_l \leq_K \bold q$
and $U_{\bold q}=U_{\bold{q}_1}+_{U_{\bold{q}_0}} U_{\bold{q}_2}$.
We shall denote this $\bold q$ by $\bold{q}_1+_{\bold{q}_0} \bold{q}_2$.

g. If $\alpha \in U_{\bold q}^+$, then $\mathbb{P}_{\bold q,\alpha}$
is a $(<\theta)$-complete forcing satisfying $*_{\theta}^{\epsilon}$
(hence $\theta^+$-cc).

h. Suppose that $\bold q \in K$ and $\underset{\sim}{\mathbb Q}$
is a $\mathbb{P}_{\bold q, \infty}$-name of a forcing notion whose
universe is from $V$, such that the conditioncs of definitions 3(d)
and 3(e) are satisfied, then there is $\bold q' \in K$ such that
$\bold q \leq_K \bold q'$, $U_{\bold q'}=U_{\bold q} \cup \{ \gamma\}$,
$U_{\bold q'} \models \alpha<\gamma$ for every $\alpha \in U_{\bold q}$
and $\underset{\sim}{\mathbb{Q}_{\bold q',\gamma}}=\underset{\sim}{\mathbb Q}$.
$\square$

\textbf{Definition 5: }The forcing notion $\mathbb P$ will be defined
as follows: 

a. The conditions of $\mathbb P$ are the elements $\bold q$ of $K\cap H(\lambda^+)$
such that $U_{\bold q} \subseteq \lambda^+$, and for every $\beta \in U_{\bold q}$,
$\underset{\sim}{\mathbb{Q}_{\beta}}$ is a name for a forcing whose
underlying set of conditions is some $X_{\beta} \subseteq \lambda^+$.

b. Given $\bold{q}_1,\bold{q}_2 \in \mathbb P$, $\mathbb P \models "\bold{q}_1 \leq \bold{q}_2"$
iff $\bold{q}_1 \leq_K \bold{q}_2$.

c. Given a generic set $G\subseteq \mathbb P$, we let $\bold{q}_G=\cup \{ \bold q : \bold q \in G\}$.

\textbf{Claim 6: }a. $\mathbb P$ is $(<\kappa)$-strategically complete.
Moreover, it's $(<\lambda^+)$-complete and $(<\theta)$-directed
closed.

b. $\Vdash_{\mathbb P} "\bold{q}_{\underset{\sim}{G}} \in K"$, hence
$\Vdash_{\mathbb P} "\mathbb{P}_{\bold q_{\underset{\sim}{G}},\infty}$
is $(<\theta)$-directed closed and $\theta^+$-cc$"$.

c. If $\delta<\lambda^+$, $cf(\delta)>\theta$ and $(\bold{q}_{\alpha} : \alpha<\delta)$
is $\leq_{\mathbb P}$-increasing, then $\bold q:= \underset{\alpha<\delta}{\cup} \bold{q}_{\alpha}$
belongs to $\mathbb P$ and $\mathbb{P}_{\bold q}=\underset{\alpha<\delta}{\cup} \mathbb{P}_{\bold{q}_{\alpha}}$.
By $\theta^+$-c.c., $\underset{\sim}{a}$ is a canonical $\mathbb{P}_{\bold q}$-name
of a member of $[\theta]^{\theta}$ iff $\underset{\sim}{a}$ is a
canonical $\mathbb{P}_{\bold{q}_{\alpha}}$-name of a member of $[\theta]^{\theta}$
for some $\alpha<\delta$.

\textbf{Proof: }The claim follows directly from the definitions. The
fact that $\Vdash_{\mathbb P} "\bold{q}_{\underset{\sim}{G}} \in K"$
follows from the general fact that if $I$ is a directed set, $\{\bold{q}_t : t\in I\} \subseteq K$
and $s\leq_I t \rightarrow \bold{q}_s \leq_K \bold{q}_t$, then $\cup \{ \bold{q}_t : t\in I\}$
is well-defined and belongs to $K$. This also shows that $\mathbb P$
is $(<\theta)$-directed closed. $\square$

We shall now define our desired model:

\textbf{Definition 7: }a. In $V^{\mathbb P}$, let $\mathbb Q= \mathbb{P}_{\bold q_{\underset{\sim}{G}},\infty}$. 

b. Let $V_2=V^{\mathbb P \star \underset{\sim}{\mathbb Q}}$.

c. Let $V_1$ be $HOD({\mathcal{P}(\theta)}^{<\kappa} \cup V)$ inside
$V_2$.

\textbf{Claim 8: }a. $V_1 \models ZF+DC_{<\kappa}$.

b. $(Ord^{<\kappa})^{V_1}=(Ord^{<\kappa})^{V_2}$, hence $\mathcal{P}(\theta)^{V_1}=\mathcal{P}(\theta)^{V_2}$.

\textbf{Proof: }We shall prove the first part of clause (b), the rest
should be clear. Clearly, $(Ord^{<\kappa})^{V_1} \subseteq (Ord^{<\kappa})^{V_2}$.
Now let $\eta \in (Ord^{\gamma})^{V_2}$ for some $\gamma<\kappa$,
then $\eta=\underset{\sim}{\eta}[G]$ for some name $\underset{\sim}{\eta}$
of a member of $Ord^{\gamma}$, where $G\subseteq \mathbb{P} \star \underset{\sim}{\mathbb Q}$
is generic. $G=G_1 \star G_2$ where $G_1 \subseteq \mathbb P$ is
generic and $G_2 \subseteq \underset{\sim}{\mathbb Q}[G_1]$ is generic.
Working in $V[G_1]$, $\underset{\sim}{\eta}/G_1$ is a $\underset{\sim}{\mathbb Q}[G_1]$-name.
As $\underset{\sim}{\mathbb Q}[G_1]$ is $\theta^+$-cc, for every
$\beta<\gamma$ there is a maximal antichain $\{ p_{\beta,i} : i<\theta \} \subseteq \underset{\sim}{\mathbb{Q}}[G_1]$
of conditions that force a value to $\underset{\sim}{\eta} /G_1 (\beta)$.
Let $\{ \zeta_{\beta,i} : i<\theta \}$ be the set corresponding values
forced by the above conditions. Let $\Gamma=\{ \underset{\sim}{p_{\beta,i}}, \underset{\sim}{\zeta_{\beta,i}} : \beta<\gamma, i<\theta\}$
be the corresponding $\mathbb P$-names for the above objects (so
we can regard them as $\mathbb P$-names for ordinals). As there are
$<\kappa$ such names and $\mathbb{P}$ is $(<\kappa)$-strategically
complete, there is a dense set of $\bold q \in \mathbb P$ that force
values to all elements of $\Gamma$. Therefore, there is some $\bold q \in \mathbb P \cap G_1$
that forces values to all elements of $\Gamma$ (and the values forced
are necessarily $\{ p_{\beta,i}, \zeta_{\beta,i} : \beta<\gamma, i<\theta\}$).
It follows that $\{ p_{\beta,i}, \zeta_{\beta,i} : \beta<\gamma, i<\theta\} \in V$.
In $V_2$, there is a function $f: \gamma \rightarrow \theta$ such
that for every $\beta<\gamma$, $\eta(\beta)=\zeta_{\beta, f(\beta)}$.
As $f\in \mathcal{P}(\theta)^{<\kappa}$ and $\{ p_{\beta,i}, \zeta_{\beta,i} : \beta<\gamma, i<\theta\} \in V$,
it follows that $\eta \in V_1$. $\square$

\textbf{Main Claim 9: }There are no $\theta$-mad families in $V_1$.

The rest of the paper will be devoted to the proof of Claim 9. 

Suppose towards contradiction that there is a $\theta$-mad family
in $V_1$, so there is some $(\bold{q}_0, \underset{\sim}{p_0}) \in \mathbb P \star \underset{\sim}{\mathbb Q}$
forcing this statement about $\underset{\sim}{\mathcal A}$ where
$\underset{\sim}{\mathcal A}$ is a canonical $\mathbb{P} \star \underset{\sim}{\mathbb Q}$-name
of a $\theta$-mad family definable using $\underset{\sim}{\eta}$,
and $\underset{\sim}{\eta}$ is a canonical $\mathbb P \star \underset{\sim}{\mathbb Q}$-name
of a parameter (so $\underset{\sim}{\eta}=((\underset{\sim}{a_{\epsilon}} : \epsilon<\underset{\sim}{\epsilon(*)}), \underset{\sim}{x})$,
where $\Vdash"\underset{\sim}{\epsilon(*)}<\kappa"$, each $\underset{\sim}{a_{\epsilon}}$
is a $\mathbb P \star \underset{\sim}{\mathbb Q}$-name of a subset
of $\theta$ and $\Vdash "\underset{\sim}{x} \in V"$). Let $G_0 \subseteq \mathbb P$
be generic over $V$ such that $\bold{q}_0 \in G_0$. In $V[G_0]$,
$\underset{\sim}{\eta}$ is a $\mathbb{P}_{\bold{q}_{G_0}, \infty}$-name,
and by increasing $\bold{q}_0$, we may assume wlog that $p_0:=\underset{\sim}{p_0}[G_0] \in \mathbb{P}_{\bold{q}_0}$,
$x=\underset{\sim}{x}[G_0] \in V$, $\epsilon(*)=\underset{\sim}{\epsilon(*)}[G_0] \in \kappa$
and that each $\underset{\sim}{a_{\epsilon}}$ ($\epsilon<\epsilon(*)$)
is a canonical $\mathbb{P}_{\bold{q}_0}$-name of a subset of $\theta$.
Given $\bold q \in \mathbb P$ above $\bold{q}_0$, let $\mathcal{A}_{\bold q}$
be the set of canonical $\mathbb{P}_{\bold q}$-names $\underset{\sim}{a}$
such that $(\bold q, \underset{\sim}{p_0}) \Vdash_{\mathbb P \times \underset{\sim}{\mathbb Q}} "\underset{\sim}{a} \in \underset{\sim}{\mathcal A}"$,
so $\bold{q}_0 \leq \bold{q}_1 \leq \bold{q}_2 \rightarrow \mathcal{A}_{\bold{q}_1} \subseteq \mathcal{A}_{\bold{q}_2}$.
Note that if $\bold{q}_0 \leq \bold{q}_1$, $\mathbb{P}_{\bold{q}_1,\infty} \models "p_0 \leq p_1"$
and $(\bold{q}_1,p_1) \Vdash "\underset{\sim}{b} \in [\theta]^{\theta}"$,
then for some $(\bold{q}_2,\underset{\sim}{a})$ we have $\bold{q}_1 \leq_{\mathbb P} \bold{q}_2$,
$\underset{\sim}{a} \in \mathcal{A}_{\bold{q}_2}$ and $(\bold{q}_2, p_0) \Vdash "\underset{\sim}{b} \cap \underset{\sim}{a} \in [\theta]^{\theta}"$.
By extending any given $\bold{q}_1 \in \mathbb P$ above $\bold{q}_0$
in this way sufficiently many times to add witnesses for madness,
and recalling Claim 6(c), we establish that the set $\{ \bold{q}_1 : \bold{q}_0 \leq_{\mathbb P} \bold{q}_1$
and $\Vdash_{\mathbb{P}_{\bold{q}_1}} "\mathcal{A}_{\bold{q}_1}$
is $\theta$-mad$"\}$ is dense in $\mathbb P$ above $\bold{q}_0$.

Now, in $V_2$, let $I=\{A \subseteq \theta: A$ is contained in a
union of $<\theta$ members of $\mathcal A \}$, then $I$ is a $\theta$-complete
ideal and $\theta \notin I$. Let $F$ be the dual filter of $I$,
then $F$ is $\theta$-complete, and as $\theta$ is supercompact
in $V_2$ (recalling that $\theta$ is Laver indestructible and that
$\mathbb P \star \underset{\sim}{\mathbb Q}$ is $(<\theta)$-directed
closed), there is a $\mathbb P \star \underset{\sim}{\mathbb Q}$-name
$\underset{\sim}{D}$ such that $(\bold{q}_0, p_0) \Vdash_{\mathbb P \star \underset{\sim}{\mathbb Q}} "\underset{\sim}{D}$
is a $\theta$-complete ultrafilter on $\theta$ that extends $F$,
and hence is disjoint to $\underset{\sim}{\mathcal A}"$. By Claim
6 and a previous observation, we may assume wlog that $\bold{q}_0 \Vdash_{\mathbb P}"\mathcal{A}_{\bold{q}_0}$
is $\theta$-mad and $\underset{\sim}{D_{\bold{q}_0}}:=\underset{\sim}{D} \cap \mathcal{P}(\theta)^{V^{\mathbb{P}_{\bold{q}_0, \infty}}}$
is a $\mathbb{P}_{\bold{q}_0,\infty}$-name of an ultrafilter on $\theta"$. 

Given an ultrafilter $U$ on $\theta$, the forcing $\mathbb{Q}_U$
is defined as follows: the conditions of $\mathbb{Q}_U$ have the
form $(u, A)$ where $u\in [\theta]^{<\theta}$ and $A\in U$. the
order is defined naturally, i.e. $(u_1, A_1) \leq (u_2, A_2)$ iff
$u_1 \subseteq u_2$, $u_2 \setminus u_1 \subseteq A_1$ and $A_2 \subseteq A_1$.

We may assume wlog that $\mathbb{P}_{\bold{q}_0, \infty}$ forces
$2^{\theta}=\lambda$, hence there is a canonical $\mathbb{P}_{\bold{q}_0, \infty}$-name
$\underset{\sim}{f}$ of a bijection from $\mathbb{Q}_{\underset{\sim}{D}_{\bold{q}_0}}$
onto $\lambda$. Let $\underset{\sim}{\mathbb Q'}$ be a name for
the forcing such that $\Vdash_{\mathbb{P}_{\bold{q}_0}} "\underset{\sim}{f}$
is an isomorphism from $\mathbb{Q}_{\underset{\sim}{D}_{\bold{q}_0}}$
onto $\underset{\sim}{\mathbb Q'}"$. Let $\underset{\sim}{B}=\underset{\sim}{B_{\underset{\sim}{D}_{\bold{q}_0}}}$
be the $\mathbb{Q}_{\underset{\sim}{D}_{\bold{q}_0}}$-name $\cup \{ u : (u,A) \in G_{\mathbb{Q}_{\underset{\sim}{D}_{\bold{q}_0}}} \}$,
so $\Vdash_{\mathbb{P}_{\bold{q}_0, \infty} \star \mathbb{Q}_{\underset{\sim}{D}_{\bold{q}_0}}} "\underset{\sim}{B} \in [\theta]^{\theta}$
is $\theta$-almost disjoint to $\mathcal{A}_{\bold{q}_0}"$. Let
$\underset{\sim}{B}'$ be the canonical $\mathbb{P}_{\bold{q}_0,\infty} \star \mathbb{Q}_{\underset{\sim}{D}_{\bold{q}_0}}$-name
for the image of $\underset{\sim}{B}$ under $\underset{\sim}{f}$. 

Now observe that there is $\bold q' \in \mathbb P$ such that $\bold{q}_0 \leq_{\mathbb P} \bold q'$,
$U_{\bold q'}=U_{\bold{q}_0} \cup \{ \gamma\}$, $\alpha<_{U_{\bold q'}} \gamma$
for every $\alpha \in U_{\bold{q}_0}$ and $\underset{\sim}{\mathbb{Q}_{\bold q', \gamma}}=\underset{\sim}{\mathbb Q'}$.
As before, there is $\bold{q}'' \in \mathbb P$ above $\bold q'$
such that $p_0 \Vdash_{\mathbb{P}_{\bold{q}'',\infty}} "\mathcal{A}_{\bold{q}''}$
is $\theta$-mad$"$. Therefore, there is some canonical $\mathbb{P}_{\bold{q}''}$-name
$\underset{\sim}{A} \in \mathcal{A}_{\bold{q}''}$ such that $p_0 \Vdash_{\mathbb{P}_{\bold{q}'',\infty}} "\underset{\sim}{A} \cap \underset{\sim}{B}' \in [\theta]^{\theta}$,
so $\underset{\sim}{A}$ has intersection of size $\theta$ with every
member of $\underset{\sim}{D}_{\bold{q}_0}$ and $\underset{\sim}{A} \notin \mathcal{A}_{\bold{q}_0}"$. 

Now let $(\bold{q}_1, \underset{\sim}{B_1}, \underset{\sim}{A_1})=(\bold q'', \underset{\sim}{B'},\underset{\sim}{A})$
and let $(\bold{q}_2, \underset{\sim}{B_2}, \underset{\sim}{A_2})$
be an isomorphic copy of $(\bold{q}_1, \underset{\sim}{B_1}, \underset{\sim}{A}_1)$
over $\bold{q}_0$ such that $U_{\bold{q}_1} \cap U_{\bold{q}_2}=U_{\bold{q}_0}$
and $\bold{q}_2 \in \mathbb P$.

\textbf{Claim 10: }Let $\bold{q}_0$, $(\bold{q}_1, \underset{\sim}{B_1}, \underset{\sim}{A_1})$
and $(\bold{q}_2, \underset{\sim}{B_2}, \underset{\sim}{A_2})$ be
as above (so $\bold{q}_0 \leq_K \bold{q}_l$ $(l=1,2)$, $U_{\bold{q}_1} \cap U_{\bold{q}_2}=U_{\bold{q}_0}$
and $\underset{l=1,2}{\wedge} \Vdash_{\mathbb{P}_{\bold{q}_l,\infty}}"\underset{\sim}{A_l} \in \underset{\sim}{\mathcal A} \setminus \mathcal{A}_{\bold{q}_0}"$)
and let $G\subseteq \mathbb{P}_{\bold{q}_0,\infty}$ be generic over
$V$, then $\Vdash_{\mathbb{P}_{\bold{q}_1,\infty} /G \times \mathbb{P}_{\bold{q}_2,\infty} /G} "\underset{\sim}{A_2} \setminus \underset{\sim}{A_1}, \underset{\sim}{A_1} \setminus \underset{\sim}{A_2} \in [\theta]^{\theta}"$.

\textbf{Proof: }We shall prove the claim for $A_2 \setminus A_1$,
the other case is similar. Suppose towards contradiction that $(p_1, p_2)$
forces that $\underset{\sim}{A_2} \setminus \underset{\sim}{A_1} \subseteq \gamma<\theta$.
For $l\in \{1, 2\}$, let $B_{l}=\{ \epsilon <\theta : p_l \nVdash_{\mathbb{P}_{\bold{q}_l},\infty /G} "\epsilon \notin \underset{\sim}{A_l}"\} \in V[G]$.
By the assumption of the claim, $B_l \in [\theta]^{\theta}$. By the
$\theta$-madness of $\underset{\sim}{\mathcal{A}_0}[G]$ in $V[G]$,
there is some $Y \in \underset{\sim}{\mathcal{A}_0}[G]$ such that
$|Y\cap B_2|=\theta$. As $p_1 \Vdash_{\mathbb{P}_{\bold{q}_1,\infty} /G} "|\underset{\sim}{A_1} \cap Y|<\theta"$,
there are $q_1$ and $\beta_1<\theta$ such that $p_1 \leq q_1 \in \mathbb{P}_{\bold{q}_1, \infty} /G$
and $q_1 \Vdash_{\mathbb{P}_{\bold{q}_1, \infty} /G} "\underset{\sim}{A_1} \cap Y \subseteq \beta_1"$.
Let $\beta_2 \in Y\cap B_2$ such that $max\{ \gamma, \beta_1 \} < \beta_2$
(recalling that $|Y\cap B_2|=\theta$). By the definition of $B_2$,
there is $q_2 \in \mathbb{P}_{\bold{q}_2, \infty} /G$ above $p_2$
that forces $"\beta_2 \in \underset{\sim}{A_2}"$. Therefore, $(p_1, p_2) \leq (q_1, q_2) \in \mathbb{P}_{\bold{q}_1, \infty} /G \times \mathbb{P}_{\bold{q}_2, \infty} /G$
and $(q_1, q_2) \Vdash_{\mathbb{P}_{\bold{q}_1, \infty} /G \times \mathbb{P}_{\bold{q}_2, \infty} /G} "\beta_2 \in \underset{\sim}{A_2} \setminus \underset{\sim}{A_1}"$,
a contradiction. It follows that $\Vdash_{\mathbb{P}_{\bold{q}_1,\infty} /G \times \mathbb{P}_{\bold{q}_2,\infty} /G} "\underset{\sim}{A_2} \setminus \underset{\sim}{A_1} \in [\theta]^{\theta}"$.
$\square$

\textbf{Claim 11: }Under the assumptions of Claim 10 (recalling that
$\Vdash_{\mathbb{P}_{\bold{q}_l, \infty}} "\underset{\sim}{A_l} \cap B \neq \emptyset$
for every $B\in \underset{\sim}{D_{\bold{q}_0}}"$ $(l=1,2)$), we
have $\Vdash_{\mathbb{P}_{\bold{q}_1,\infty} /G \times \mathbb{P}_{\bold{q}_2,\infty} /G} "\underset{\sim}{A_1} \cap \underset{\sim}{A_2} \in [\theta]^{\theta}"$.

\textbf{Proof: }Assume towards contradiction that $(p_1, p_2) \in \mathbb{P}_{\bold{q}_1,\infty} /G \times \mathbb{P}_{\bold{q}_2,\infty} /G$
forces that $\underset{\sim}{A_1} \cap \underset{\sim}{A_2} \subseteq \gamma$
for some $\gamma<\theta$. It's forced by $(p_1, p_2)$ that $\underset{\sim}{A_l} \subseteq B_l$
$(l=1,2)$ where $B_l$ is as in the proof of the previous claim,
hence it's forced by $(p_1,p_2)$ that each $B_l$ intersects each
member of $\underset{\sim}{D_{\bold{q}_0}}$. As $B_1, B_2 \in V[G]$,
it follows that $B_1,B_2 \in \underset{\sim}{D_{\bold{q}_0}}[G]$.
Therefore, there is some $\beta \in (B_1 \cap B_2) \setminus \gamma$,
hence there is $q_l \in \mathbb{P}_{\bold{q}_l, \infty} /G$ above
$p_l$ that forces $"\beta \in \underset{\sim}{A_l}"$ $(l=1,2)$.
It follows that $(p_1, p_2) \leq (q_1, q_2) \in \mathbb{P}_{\bold{q}_1,\infty} /G \times \mathbb{P}_{\bold{q}_2,\infty} /G$
and $(q_1, q_2) \Vdash_{\mathbb{P}_{\bold{q}_1,\infty} /G \times \mathbb{P}_{\bold{q}_2,\infty} /G} "\beta \in \underset{\sim}{A_1} \cap \underset{\sim}{A_2}"$,
contradicting the choice of $\gamma$ and $(p_1, p_2)$. It follows
that $\Vdash_{\mathbb{P}_{\bold{q}_1,\infty} /G \times \mathbb{P}_{\bold{q}_2,\infty} /G} "\underset{\sim}{A_1} \cap \underset{\sim}{A_2} \in [\theta]^{\theta}"$.
$\square$

Now given $\bold{q}_0$, $(\bold{q}_1, \underset{\sim}{B_1}, \underset{\sim}{A_1})$
and $(\bold{q}_2, \underset{\sim}{B_2}, \underset{\sim}{A_2})$ as
above, let $\bold{q}_3=\bold{q}_1+_{\bold{q}_0} \bold{q}_2$. Then
$\bold{q}_3 \in \mathbb P$, $\bold{q}_1, \bold{q}_2 \leq_K \bold{q}_3$,
and by claims 10 and 11, we get a contradiction. This completes the
proof of Main Claim 9 and hence of Theorem 2. $\square$

$\\$

We conclude with the following natural question:

\textbf{Question: }What's the consistency strength of $ZF+DC_{\theta}+"$there
are no $\theta$-mad families$"$ for some $\theta>\aleph_0$?

$\\$

\textbf{\large References}{\large \par}

{[}HwSh:1090{]} Haim Horowitz and Saharon Shelah, Can you take Toernquist's
inaccessible away? arXiv:1605.02419

{[}HwSh:1093{]} Haim Horowitz and Saharon Shelah, Transcendence bases,
well-orderings of the reals and the axiom of choice, arXiv:1901.01508

{[}HwSh:1113{]} Haim Horowitz and Saharon Shelah, Madness and regularity
properties, arXiv:1704.08327

{[}HwSh:1145{]} Haim Horowitz and Saharon Shelah, $\kappa$-Madness
and Definability, arXiv:1805.07048

{[}Ma{]} A. R. D. Mathias, Happy families, Ann. Math. Logic \textbf{12
}(1977), no. 1, 59-111

{[}Sh:218{]} Saharon Shelah, On measure and category, Israel J. Math.
52 (1985) 110-114

{[}To{]} Asger Toernquist, Definability and almost disjoint families,
Advances in Mathematics 330, 61-73, 2018

{[}Vi{]} Matteo Viale, Category forcings, $MM^{+++}$, and generic
absoluteness for the theory of strong forcing axioms, J. Amer. Math.
Soc. 29 (2016), no. 3, 675-728

$\\$

(Haim Horowitz) Department of Mathematics

University of Toronto

Bahen Centre, 40 St. George St., Room 6290

Toronto, Ontario, Canada M5S 2E4

E-mail address: haim@math.toronto.edu

$\\$

(Saharon Shelah) Einstein Institute of Mathematics

Edmond J. Safra Campus,

The Hebrew University of Jerusalem.

Givat Ram, Jerusalem 91904, Israel.

Department of Mathematics

Hill Center - Busch Campus, 

Rutgers, The State University of New Jersey.

110 Frelinghuysen Road, Piscataway, NJ 08854-8019 USA

E-mail address: shelah@math.huji.ac.il
\end{document}